\documentclass[a4paper,12pt]{article}
\usepackage{fancyhdr,graphicx}
\usepackage{amsfonts}
\usepackage{amssymb}
\usepackage{amsthm}
\usepackage{newlfont}
\usepackage{amsmath}
\usepackage[top=2.0cm,bottom=3.0cm,left=2.5cm,right=2.5cm]{geometry}
\usepackage{multirow}
\usepackage{array}
\usepackage{longtable}
\usepackage{bm}
\usepackage{latexsym}
\usepackage{amsmath}
\usepackage{graphicx}
\usepackage{amssymb}
\usepackage{mathrsfs}
\usepackage{setspace}
\allowdisplaybreaks \numberwithin{equation}{section}

\begin{document}

\title{{O-Fibonacci $(p,r)$-cube as Cartesian products
\thanks{This work is supported by the Natural Science Foundation of ShanDong Province (ZR2018MA010),
the National Natural Science Foundation of China (Grant no. 11861032).}}}

\author{Jianxin Wei$^{a,}$\footnote{Corresponding author.}, Guangfu Wang$^{b}$\\
\scriptsize{$^{a}$
School of Mathematics and Statistics Science, Ludong University, Yantai, Shandong, 264025, PR China}\\
  \scriptsize{$^{b}$
School of Science, East China Jiaotong University, Nanchang, Jiangxi, 330013, PR China}\\
\scriptsize{E-mail addresses: wjx0426@163.com(J. Wei),  wgfmath@126.com(G. Wang)}}

\date{}

\maketitle

\begin{abstract}
Let $p ,r $ and $n $ be positive integers.
Then the O-Fibonacci $(p,r)$-cube $O\Gamma^{(p,r)}_{n}$ is the subgraph of $Q_{n}$ induced on the binary words in which there is at least $p-1$ zeros between any two $1$s and there is at most $r$ consecutive $10^{p-1}$.
These cubes include a wide range of cubes as their special cases,
such as hypercubes, Fibonacci cubes, and postal networks.
In this note it is proved that $O\Gamma^{(p,r)}_{n}$ is a non-trivial Cartesian product if and only if $p=1$ and $r\geq n\geq2$.
\end{abstract}

\newcommand{\trou}{\vspace{1.5 mm}}
\newcommand{\noi}{\noindent}
\newcommand{\ol}{\overline}

\section{Introduction}
A hypercube $Q_{n}$ can be defined as the graph whose vertex set consists of all binary words of length $n$,
where two vertices are adjacent if and only if they differ in precisely one coordinate.
The cube $Q_{3}$ is shown in Fig. 1$(a)$.
For $n\geq1$, 
a Fibonacci cube $\Gamma_{n}$ is the graph obtained from $Q_{n}$ by removing all vertices that contain no two consecutive $1$s  \cite{Hsu}. 
The cube $\Gamma_{5}$ is shown in Fig. 1$(b)$.
For more results on application and structure of  Fibonacci cubes, see \cite{K1} for a survey.

Based on the Fibonacci $(p,r)$-numbering system,
Egiazarian and Astola\cite{Egiazarian} defined the $O$-Fibonacci $(p,r)$-cubes as follows.
For $n\geq1$,
$\alpha=a_{1}a_{2}\ldots a_{n}$ is called a $O$-Fibonacci $(p,r)$-word if the following hold:
(1) if $b_{i}=1$ then $b_{i+1}=\ldots=b_{i+(p-1)}=0$,
i.e. there is at least $(p-1)$ 0s between two 1s (which is called `consecutive' 1s);
(2) there are no more than $r$ `consecutive' 1s in $\alpha$,
i.e. there are at most $r$ consecutive $10^{p-1}$ in $\alpha$.

\trou \noi {\bf Definition 1.1\cite{Egiazarian}.} \emph{For positive integers $p$, $r$ and $n$,
let $O\mathcal{F}_{n}^{(p,r)}$ be the set of all the $O$-Fibonacci $(p,r)$-word of length $n$.
Then the $O$-Fibonacci $(p,r)$-cube $O\Gamma^{(p,r)}_{n}$ is the graph defined on the vertex set $O\mathcal{F}_{n}^{(p,r)}$,
and two vertices being adjacent if they differ exactly in one coordinate.}

\trou \noi {\bf Definition 1.1$'$.}
\emph{Let $p$, $r$ and $n$ be any positive integers.
Then the $O$-Fibonacci $(p,r)$-cube $O\Gamma^{(p,r)}_{n}$ is a subgraph of $Q_{n}$ induced on vertices}

$V(O\Gamma_{n}^{(p,r)})=0V(O\Gamma_{n-1}^{(p,r)})\cup 10^{p-1}0V(O\Gamma_{n-p-1}^{(p,r)})\cup\ldots \cup (10^{p-1})^{r}0V(O\Gamma_{n-pr-1}^{(p,r)})$

\noindent
\emph{with the initial conditions $V(O\Gamma_{i}^{(p,r)})=\emptyset$ for $i<0$,
$V(O\Gamma_{i}^{(p,r)})=\{\lambda\}$ for $i=0$}.

The cubes $O\Gamma^{(1,3)}_{3}$,
$O\Gamma^{(1,1)}_{5}$, $O\Gamma^{(2,1)}_{6}$,
$O\Gamma^{(2,2)}_{5}$ and $O\Gamma^{(2,2)}_{4}$ are shown in Fig. 2$(a)$,$(b)$,$(c)$,$(d)$ and $(f)$,
respectively.

\begin{figure}[h]
\begin{center}
\includegraphics[scale=0.80]{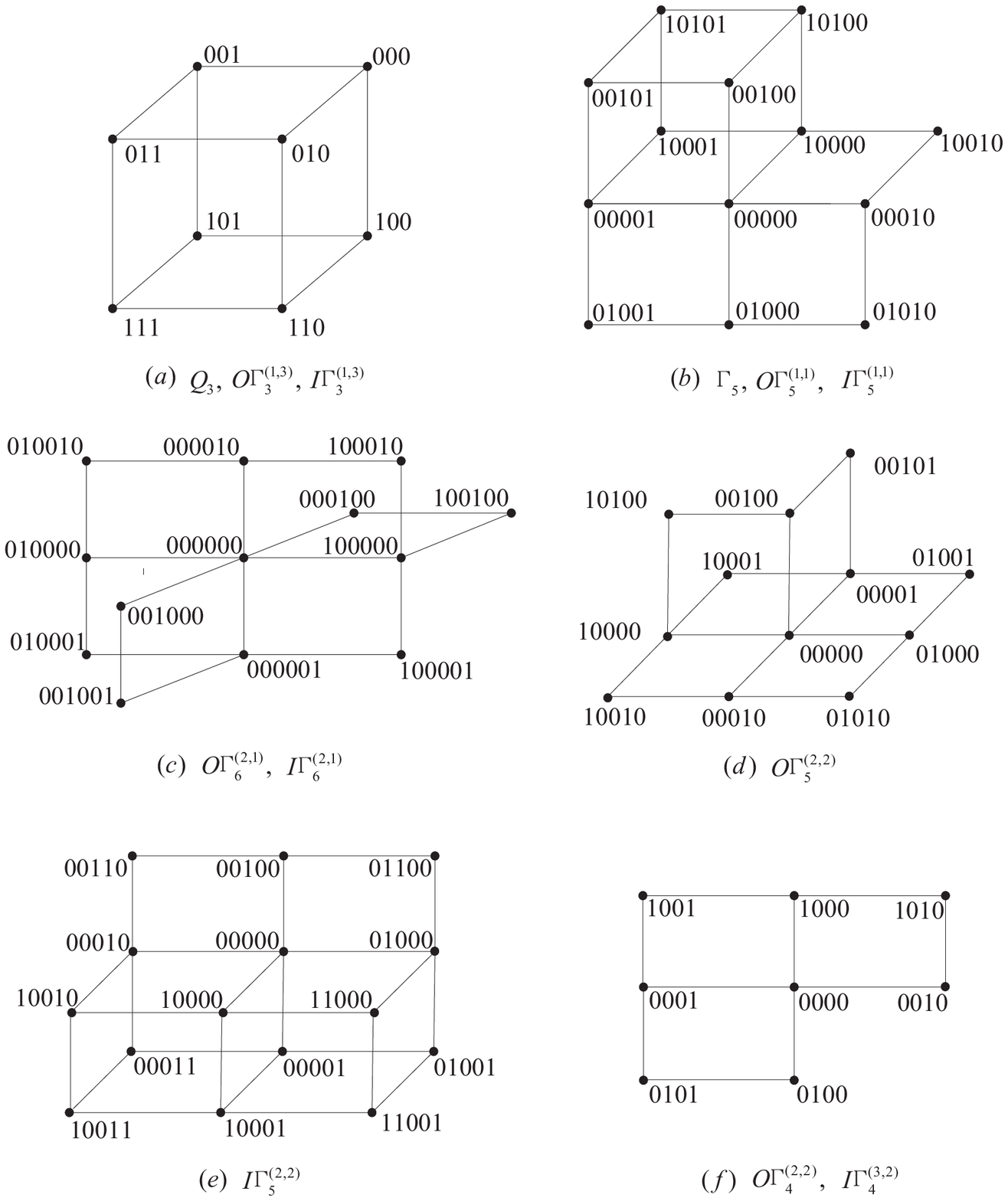}
\\{Fig. 1. Some examples of $Q_{n}$, $\Gamma_{n}$, $O\Gamma^{(p,r)}_{n}$ and $I\Gamma^{(p,r)}_{n}$.}
\end{center}
\end{figure}

In papers \cite{K1,KR,OZ,OZY,JZ},
another ``Fibonacci $(p,r)$-cube'' is studied from several different aspects.
It was shown that this cube is a new topological structure differing from the existing Fibonacci like-cubes \cite{JYW}.

Let $p$, $r$ and $n$ be positive integers.
Then a $I$-Fibonacci $(p,r)$-word of length $n$ is a word of length $n$
in which there are at most $r$ consecutive 1s and at least $p$ element $0$s between two sub-words composed of (at most $r$) consecutive 1s.

\trou \noi {\bf Definition 1.2\cite{OZY}.}
\emph{Let $I\mathcal{F}_{n}^{(p,  r)}$ denote the set of $I$-Fibonacci $(p,r)$-codes of length $n$.
Then $I\Gamma_{n}^{(p,r)}=(V,E)$ is the graph such that $V=I\mathcal{F}_{n}^{(p,r)}$ and two vertices are adjacent if they differ in exactly one coordinate.}

For examples,
the cubes $I\Gamma_{3}^{(1,3)}$, $I\Gamma_{5}^{(1,1)}$, $I\Gamma_{6}^{(2,1)}$, $I\Gamma_{5}^{(2,2)}$ and $I\Gamma_{4}^{(3,2)}$ are shown in Fig. 1$(a)$, $(b)$, $(c)$, $(e)$ and $(f)$, respectively.
Although $O\Gamma_{n}^{(p,r)}$ and $I\Gamma_{n}^{(p,r)}$ are different cubes in general,
there are also exist the same cubes defined by different $p,r$ and $n$, see Fig. 1$(f)$.

\trou \noi {\bf Property 1.3\cite{JYW}.}
\emph{Let $r\geq1$, $p\geq1$ and $n\geq1$.
Then $O\Gamma_{n}^{(p,r)}\cong I\Gamma_{n}^{(p,r)}$ if and only if $p=1$ or $r=1$.}

This note was primarily motivated with the study of \cite{KR} which give the characterization of $I\Gamma_{n}^{(p,r)}$ that are Cartesian products.
Here we also study this question for $O\Gamma_{n}^{(p,r)}$,
which is also a problem posed in \cite{JYW}.
We get the following result applying the similar method in \cite{KR,OZY}:

\trou \noi {\bf Theorem 1.4.}
\emph{Let $p\geq1$, $r\geq1$ and $n\geq1$.
Then $O\Gamma_{n}^{(p,r)}$ is a non-trivial Cartesian product graph if and only if $p=1$ and $r\geq n\geq2$.}

The rest of paper is organized as follows.
In the next section some preliminary definitions and results are given.
In the last section, Theorem 1.3 is proved.

\section{Preliminaries}

Let $u$ and $v$ be any two binary words.
Then $uv$ denotes its its concatenation.
With $u^{n}$ we mean the concatenation of $n$ copies of $u$.
For example,
$1^{n}$ is the binary word of length $n$,
and $u^{0}$ is the empty word $\lambda$.
For a word $\alpha=a_{1}a_{2}\ldots a_{n}$,
$w(\alpha)=\sum^{n}_{i=1}a_{i}$ is called the weight of $\alpha$,
in other words,
$w(\alpha)$ is the number of $1$s in $\alpha$.
A word with $1$ in coordinate $i$ and $0$ elsewhere,
denoted with $e^{i}$,
is called the $i$th unit word.

The distance $d_{G}(u, v)$ between vertices $u$ and $v$ of a connected graph $G$ is the length of a shortest $u, v$-path.
Sometimes we simply write $d(u, v)$ instead of $d_{G}(u, v)$.
The Hamming distance $H(u, v)$ between binary words $u$ and $v$ (of the same length) is the number of coordinates in which
they differ. It is well known that $d_{Q_{n}} (u, v) = H(u, v)$ holds for any $u, v \in V(Q_{n})$.

The Cartesian product of $G$ and $H$ is a graph,
denoted as $G\Box H$,
whose vertex set is $V(G)\times V(H)$.
Two vertices $(g, h)$ and $(g', h')$ are adjacent precisely if $g = g'$ and $hh'\in E(H)$,
or $h = h'$ and $gg'\in E(G)$.
The graphs $G$ and $H$ are called factors of the product $G\Box H$.
If $G= G_{1}\Box G_{2}\Box \ldots \Box G_{k}$,
then we say $G_{1}\Box G_{2}\Box \ldots \Box G_{k}$ is a product representation of $G$.
A graph is prime with respect to the Cartesian product if it is nontrivial and cannot be
represented as the product of two nontrivial graphs.
For more information on the Cartesian product of graphs see \cite{HIK}.
It is well known that the hypercube $Q_{n}$ is the Cartesian product of $n$ copies of $K_{2}$.

Recall that two edges $e=ab$ and $f=xy$ in a graph are in relation $\Theta$,
in symbols $e\Theta f$, if $d(a,x)+ d(b,y)\neq d(a,y)+d(b,x)$.
The relation $\Theta$ is reflexive and symmetric, but need not be transitive.
For a connected graph $G= G_{1}\Box G_{2}\Box \ldots \Box G_{k}$,
the product color $c(uv)$ of an edge $uv$ of $G$ is $i$ if $u$ and $v$ differ in coordinate $i$,
and two edges $e$ and $f$ of $G$ are called in product relation if $c(e)=c(f)$.
The product relation is transitive, reflexive, and symmetric.
For any index $i$ $(i=1,2,\ldots,k)$,
there is a projection map
$pi :G_{1}\Box G_{2}\Box \ldots \Box G_{k}\rightarrow G_{i}$ defined as $p_{i}(x_{1}, x_{2}, \ldots , x_{k}) = x_{i}$.
We call $x_{i}$ the $i$th coordinate of the vertex $(x_{1}, x_{2},\ldots , x_{k})$.
The following result can be found in Lemma 13.5(i) of \cite{HIK}.

\trou \noi {\bf Lemma 2.1\cite{HIK}.}
\emph{Suppose $G= G_{1}\Box G_{2}\Box \ldots \Box G_{k}$ and $e,f\in E(G).$
If $c(e)=c(f)=i$ and $p_{i}(e)=p_{i}(f)$,
then $e\Theta f$.}

Two edges $e=uv$ and $f=uw$ are called in relation $\tau$,
in symbols $e\tau f$,
if $u$ is the unique common neighbor of $v$ and $w$.
For a relation $R$,
let $R^{*}$ be the transitive closure of R.
Then the following fundamental theorem give a tool to prove the main result of this note.

\trou \noi {\bf Theorem 2.2\cite{TF,HIK}.}
\emph{If $G$ is a connected graph, then $(\Theta \cup \tau)^{*}$ is a product relation}.

This result immediately implies that a graph is prime if and only if the
relation$(\Theta \cup \tau)^{*}$ has a single equivalence class.

Finally, some notations are introduced.
For a graph $G$£¬
let $G(X)$ denote the subgraph of $G$ induced on $X\subseteq V(G)$.
For two subsets $X_{1},X_{2}\subseteq V(G)$,
let the notation $[X_{1},X_{2}]$ denote the set of edges of $G$ having one end vertex in $X_{1}$ and the other in $X_{2}$.
The set $\{1,2,\ldots,n\}$ is denoted with $[n]$, and
the disjoint union of sets with $\uplus$.
It is written as $G\cong H$ to denote that graphs $G$ and $H$ are isomorphic graphs.

\section{Proof of Theorem 1.4}

Let $p,r,n$ be positive integers and $O\Gamma^{(p,r)}_{n}$ the $O$-Fibonacci $(p,r)$-cube.
Then the following result holds by Definition 1.1.

\trou \noi {\bf Lemma 3.1.}
\emph{Let $p\geq1,r\geq$, and $n\geq1$. Then we have}:

$(i)$ \emph{if $p=1$ and $n\leq r$, then $O\Gamma^{(p,r)}_{n}\cong Q_{n}$; and}

$(ii)$ \emph{if $r=1$ and $n\leq p+1$, or $r\geq2$ and $n\leq p$,
then $O\Gamma^{(p,r)}_{n}\cong K_{1,n}$.}

\trou \noi {\bf Proof.} If $p=1$ and $n\leq r$,
then $V(O\Gamma^{(p,r)}_{n})=V(Q_{n})$ and so $(i)$ holds.
If $r=1$ and $n\leq p+1$, or $r\geq2$ and $n\leq p$,
then $V(O\Gamma^{(p,r)}_{n})=\{u:$ $u$ is a word of length $n$ and $w(u)\leq1\}.$
Hence, $O\Gamma^{(p,r)}_{n}\cong K_{1,n}$.
$\Box$

For $i\in [n]$,
Let $E^{(p,r)}_{i}$ be the set of edges of $O\Gamma^{(p,r)}_{n}$ in whose endpoints differ in exactly coordinate $i$,
i.e.,

$E^{(p,r)}_{i}$=$\{xy\in E(O\Gamma^{(p,r)}_{n}): x, y \in O\mathcal{F}_{n}^{(p,r)}$ differ in precisely the $i$th coordinate$\}$, and set

$V^{0}_{i} = \{u : u$ is an end vertex of an edge from $E^{(p,r)}_{i}$ with $u_{i} = 0\}$,

$V^{1}_{i} = \{u : u$ is an end vertex of an edge from $E^{(p,r)}_{i}$ with $u_{i} = 1\}$.

Then $V^{0}_{i}, V^{1}_{i} \in O\mathcal{F}_{n}^{(p,r)}$ and $E^{(p,r)}_{i}=[V^{0}_{i},V^{1}_{i}]$.
For these sets the following lemma holds.

\trou \noi {\bf Lemma 3.2.}
\emph{Let $i\in [n]$. Then}

$(i)$ \emph{$O\Gamma^{(p,r)}_{n}[V^{0}_{i}]\cong O\Gamma^{(p,r)}_{n}[V^{1}_{i}]$,
$O\Gamma^{(p,r)}_{n}[V^{0}_{i} \uplus V^{1}_{i}]\cong O\Gamma^{(p,r)}_{n}[V^{0}_{i}]\Box K_{2}$ and},

$(ii)$ \emph{$O\Gamma^{(p,r)}_{n}[V^{0}_{i}]$ and $O\Gamma^{(p,r)}_{n}[V^{1}_{i}]$ are connected subgraphs of $O\Gamma^{(p,r)}_{n}$}.

\trou \noi {\bf Proof.} Since the vertices of $V^{0}_{i}$ and $V^{1}_{i}$ differ in precisely the $i$th coordinate,
$O\Gamma^{(p,r)}_{n}[V^{0}_{i}]\cong O\Gamma^{(p,r)}_{n}[V^{1}_{i}]$ holds obviously.
With the fact $E^{(p,r)}_{i}=[V^{0}_{i},V^{1}_{i}]$,
we know $E(O\Gamma^{(p,r)}_{n}[V^{0}_{i} \uplus V^{1}_{i}])$
$=E(O\Gamma^{(p,r)}_{n}[V^{0}_{i}])\uplus E^{(p,r)}_{i}\uplus E(O\Gamma^{(p,r)}_{n}[V^{1}_{i}])$.
So $O\Gamma^{(p,r)}_{n}[V^{0}_{i} \uplus V^{1}_{i}]\cong O\Gamma^{(p,r)}_{n}[V^{0}_{i}]\Box K_{2}$.
It complete the proof of $(i)$.

As $O\Gamma^{(p,r)}_{n}[V^{0}_{i}]\cong O\Gamma^{(p,r)}_{n}[V^{1}_{i}]$,
we only need to show that $O\Gamma^{(p,r)}_{n}[V^{0}_{i}]$ is connected for any $i\in [n]$.
It is clear that $0^{n}\in V^{0}_{i}$.
Let $\alpha=a_{1}a_{2}\ldots a_{n}\neq0$ be any vertex of $V^{0}_{i}$ such that $w(\alpha)=k$ and $a_{t_{1}}=a_{t_{2}}\ldots a_{t_{k}}=1$.
By Definition 1.1$'$,
all the words $\alpha_{1}=\alpha+e^{t_{1}}$,
$\alpha_{2}=\alpha+e^{t_{1}}+e^{t_{2}}$,
$\ldots$, and
$\alpha_{k}=0^{n}=\alpha+e^{t_{1}}+e^{t_{2}}+\ldots+e^{t_{k}}$ are vertices of $V^{0}_{i}$.
So there exist a path between $\alpha$ and $0^{n}$:
$\alpha$-$\alpha_{1}$-$\alpha_{2}$-\ldots-$\alpha_{k}$ in $O\Gamma^{(p,r)}_{n}[V^{0}_{i}]$.
This means that any vertex of $V^{0}_{i}$ is connected with $0^{n}$.
So $O\Gamma^{(p,r)}_{n}[V^{0}_{i}]$ is connected.
$\Box$

\trou \noi {\bf Corollary 3.3.} \emph{Let $i\in [n]$ and $e, f \in E^{(p,r)}_{i}$. Then $e\Theta f$.}

\trou \noi {\bf Proof.}
Since $e, f \in E^{(p,r)}_{i}$, $c(e)=c(f)=i$ obviously.
Further by Lemma 3.2,
we know that $O\Gamma^{(p,r)}_{n}[V^{0}_{i} \uplus V^{1}_{i}]\cong O\Gamma^{(p,r)}_{n}[V^{0}_{i}]\Box K_{2}$.
By the facts $E^{(p,r)}_{i}=[V^{0}_{i},V^{1}_{i}]$ and $O\Gamma^{(p,r)}_{n}[V^{0}_{i}]\cong O\Gamma^{(p,r)}_{n}[V^{1}_{i}]$, $p_{i}(e)=p_{i}(f)=K_{2}$ holds.
Hence,
$e\Theta f$ by Lemma 2.1.
$\Box$

\trou \noi {\bf Proof of Theorem 1.3.}
If $n=1$,
then $O\Gamma^{(p,r)}_{n}\cong K_{2}$ and so it is prime.
Hence we can assume in the rest that $n\geq 2$.

By Theorem 2.2,
it suffices to show that the relation $(\Theta \cup \tau)^{*}$ consists of a single equivalence class to prove that $O\Gamma^{(p,r)}_{n}$ is prime.
For any $i\in [n]$,
all the edges of $E^{(p,r)}_{i}$ are in the same $\Theta^{*}$-class by Corollary 3.3.
The binary relation $\sim$ on the set $[n]$ also can be defined as follows\cite{KR}:
for $i,j\in [n]$ we call $i\sim j$,
if there exist edges $e\in E^{(p,r)}_{i}$ and $f\in E^{(p,r)}_{j}$ such that $e\tau f $.
Then it follows that $O\Gamma^{(p,r)}_{n}$ is a prime graph as soon as $\sim^{\ast}= [n]\times[n]$.
We distinguish two cases.

Case 1: $p=1$.

By Lemma 3.1, $O\Gamma^{(1,r)}_{n}\cong Q_{n}$ for $n\leq r$.
If $n>r$,
then $O\Gamma_{n}^{(p,r)}\cong I\Gamma_{n}^{(p,r)}$ by Property 1.3.
For $I\Gamma_{n}^{(p,r)}$ such that $n>r$,
it has shown that $\sim^{\ast}=[n]\times[n]$ in paper \cite{KR}.
Hence, $O\Gamma^{(1,r)}_{n}\cong Q_{n}$ is prime if and only if $r<n$.

Case 2: $p\geq2$.

We consider the words $e^{1}=100^{n-2}\in V^{1}_{1}$, $0^{n}=000^{n-2}$ and $e^{2}=010^{n-2}\in V^{1}_{2}$.
It is clear that $0^{n}$ is a common neighbor of $e^{1}$ and $e^{2}$,
and $e^{1}+e^{2}=110^{n}$ is another possible common neighbor of them.
Since $p\geq2$, we know that $e^{1}+e^{2}$ is not a vertex of $O\Gamma^{(p,r)}_{n}$ by Definitions 1.1 or 1.1$'$.
Hence, we know that $1\sim2$.
In general, 
we consider the words $e^{i}=100^{n-2}\in V^{1}_{1}$, $0^{n}=000^{n-2}$ and $e^{i+1}=010^{n-2}\in V^{1}_{2}$ for $i\in[n-1]$.
Then we get $i\sim i+1$,
and so $1\sim2$, $\ldots$, $n-1\sim n$.
Hence we conclude that $\sim^{\ast}=[n]\times[n]$. 
This means that for any $p\geq2$,
the cube $O\Gamma^{(p,r)}_{n}$ is prime.
$\Box$

\end{document}